# Convex Imprecise Previsions: Basic Issues and Applications

R. PELESSONI
*University of Trieste, Italy*

P. VICIG
*University of Trieste, Italy*

**Abstract**

In this paper we study two classes of imprecise previsions, which we termed convex and centered convex previsions, in the framework of Walley's theory of imprecise previsions. We show that convex previsions are related with a concept of convex natural estension, which is useful in correcting a large class of inconsistent imprecise probability assessments. This class is characterised by a condition of avoiding unbounded sure loss. Convexity further provides a conceptual framework for some uncertainty models and devices, like unnormalised supremum preserving functions. Centered convex previsions are intermediate between coherent previsions and previsions avoiding sure loss, and their not requiring positive homogeneity is a relevant feature for potential applications. Finally, we show how these concepts can be applied in (financial) risk measurement.

## 1 Introduction

Imprecise probability theory is developed by P. Walley in [14] in terms of two major classes of (unconditional) imprecise previsions, relying upon reasonable consistency requirements: *avoiding sure loss* and *coherent* previsions. The condition of avoiding sure loss is less restrictive than coherence but is often too weak.

Coherent imprecise previsions have been studied more extensively, while imprecise previsions that avoid sure loss received less attention, and it is an interesting problem to state whether some special class of previsions avoiding sure loss can be identified, which is such that

(a) its properties are not too far from those of coherent previsions;

(b) it gives further insight into the theory of imprecise previsions or generalises some of its basic aspects;

(c) it may express beliefs which do not match with coherence but which are useful in formalising and dependably modelling certain kinds of problems.

The main aim of this paper is to discuss the properties and some applications of two classes of imprecise previsions, which we termed convex and centered convex previsions and which let us provide some answers to points (a), (b), (c). The paper partly summarises and complements [12], where proofs may be found for those results which are stated without proof here.

After recalling some basic notions in Section 2, we study the larger class of convex lower previsions in Section 3.1. Although our conclusion is that convexity is an unsatisfactory consistency requirement – for instance, convex previsions do not necessarily avoid sure loss – it is however important as far as (b) is concerned. That is seen in Section 3.2, where a notion of convex natural extension is discussed which formally parallels the basic concept of natural extension in [14]. We characterise lower previsions whose convex natural extension is finite as those complying with the (mild) requirement of avoiding unbounded sure loss. In this case the convex natural extension indicates a canonical (least-committal) way of correcting them into a convex assessment. As discussed in Section 3.2.1, it is then easy to make a further correction to achieve the stronger (and more satisfactory) property of centered convexity.

Centered convex previsions are discussed in Section 3.3, together with generalisations of the important envelope theorem. Centered convex lower previsions are a special class of previsions avoiding sure loss, retaining several properties of coherent imprecise previsions, and hence they appear to fulfil requirement (a).

Section 4 gives some answers to point (c). Here convex previsions provide a conceptual framework for certain kinds of uncertainty models, as shown in Examples 1 (overly prudential assessments) and 2 (supremum preserving functions). These models are sometimes employed in practice, although they cannot usually be regarded as satisfactory. Centered convex previsions do not require the positive homogeneity condition $\underline{P}(\lambda X) = \lambda \underline{P}(X)$, $\forall \lambda > 0$, and hence seem appropriate to capture risk aversion. In Section 4.1 we focus in particular on risk measurement problems, showing that the results in Section 3 may be used to define convex risk measures (centered or not) for an arbitrary set of random variables $\mathcal{D}$. In particular, the definition of convex risk measure coincides, when $\mathcal{D}$ is a linear space, with the concept of convex risk measure recently introduced in the literature to consider liquidity risks [4, 5, 7]. It appears here that results from the risk measurement area can profitably contribute to the development of imprecise probability theory and viceversa. Section 5 concludes the paper.

## 2 Preliminaries

Unless otherwise specified, in the sequel we shall denote with $\mathcal{D}$ an *arbitrary* set of bounded random variables (or gambles, in the notation of [14]) and with $\mathcal{L}$ ($\supset \mathcal{D}$) the set of all bounded random variables (on a possibility space). A *lower prevision* $\underline{P}$ (an *upper prevision* $\overline{P}$, a *prevision* $P$) on $\mathcal{D}$ is a real-valued function with domain $\mathcal{D}$. In particular, if $\mathcal{D}$ contains only indicator functions of events, $\underline{P}$ ($\overline{P}$, $P$) is termed lower probability (upper probability, probability).

Lower (and upper) previsions should satisfy some consistency requirements: the condition of *avoiding sure loss* and the stronger *coherence* condition [14].

**Definition 1** *$\underline{P}: \mathcal{D} \to \mathbb{R}$ is a lower prevision on $\mathcal{D}$ that* avoids sure loss *iff, for all $n \in \mathbf{N}^+$, $\forall X_1, \ldots, X_n \in \mathcal{D}$, $\forall s_1, \ldots, s_n$ real and non-negative, defining $\underline{G} = \sum_{i=1}^n s_i(X_i - \underline{P}(X_i))$, $\sup \underline{G} \geq 0$.*

**Definition 2** *$\underline{P}: \mathcal{D} \to \mathbb{R}$ is a* coherent lower prevision *on $\mathcal{D}$ if and only if, for all $n \in \mathbf{N}^+$, $\forall X_0, X_1, \ldots, X_n \in \mathcal{D}$, $\forall s_0, s_1, \ldots, s_n$ real and non-negative, defining $\underline{G} = \sum_{i=1}^n s_i(X_i - \underline{P}(X_i)) - s_0(X_0 - \underline{P}(X_0))$, $\sup \underline{G} \geq 0$.*

The condition of avoiding sure loss is too weak under many respects: for instance, it does not require that $\underline{P}(X) \geq \inf X$, nor does it impose monotonicity. On the other hand, it is simpler to assess and to check than coherence.

Behaviourally, a lower prevision assessment $\underline{P}(X)$ may be viewed as a supremum buying price for $X$ [14], and $s(X - \underline{P}(X))$ represents an *elementary gain* from a bet on $X$, with stake $s$. We shall say that the bet is *in favour* of $X$ if $s \geq 0$, whilst $-s(X - \underline{P}(X))$ ($s \geq 0$) is an elementary gain from a bet *against* $X$. Definitions 1 and 2 both require that no admissible linear combination $\underline{G}$ of elementary gains originates a sure loss bounded away from zero. The difference is that the concept of avoiding sure loss considers only bets in favour of the $X_i$, while coherence considers also (at most) one bet against a random variable in $\mathcal{D}$.

We recall the following properties of coherent lower previsions, which hold whenever the random variables involved are in $\mathcal{D}$:

(a) $\underline{P}(\lambda X) = \lambda \underline{P}(X)$, $\forall \lambda > 0$ (positive homogeneity)
(b) $\inf X \leq \underline{P}(X) \leq \sup X$ (internality)
(c) $\underline{P}(X + Y) \geq \underline{P}(X) + \underline{P}(Y)$ (superlinearity).

*Coherent precise* previsions may be defined by modifying Definition 2 to allow $n \geq 0$ bets in favour of and $m \geq 0$ bets against random variables in $\mathcal{D}$ ($m, n \in \mathbf{N}$). A coherent precise prevision $P$ is necessarily *linear* and *homogeneous*: $P(aX + bY) = aP(X) + bP(Y)$, $\forall a, b \in \mathbb{R}$. In particular $P(0) = 0$.

Coherent lower previsions may be characterised using precise previsions [14]:

**Theorem 1** (Lower envelope theorem) *A lower prevision $\underline{P}$ on $\mathcal{D}$ is coherent iff $\underline{P}$ is the lower envelope of some set $\mathcal{M}$ of coherent precise previsions on $\mathcal{D}$, i.e. iff*

$$\underline{P}(X) = \inf_{P \in \mathcal{M}} \{P(X)\}, \forall X \in \mathcal{D} \text{ (inf is attained)}.$$

*Upper* and *lower* previsions are customarily related by the *conjugacy* relation $\overline{P}(X) = -\underline{P}(-X)$. An upper prevision $\overline{P}(X)$ may be viewed as an infimum selling price for $X$ and an *elementary gain* from a bet concerning $X$ is written as $s(\overline{P}(X) - X)$. The definitions of coherence and of the condition of avoiding sure loss are modified accordingly.

## 3 Convex Lower Previsions

### 3.1 Convex Previsions

**Definition 3** $\underline{P} : \mathcal{D} \to \mathbb{R}$ *is a* convex lower prevision *on $\mathcal{D}$ iff, for all $n \in \mathbf{N}^+$, $\forall X_0, X_1, \ldots, X_n \in \mathcal{D}$, $\forall s_1, \ldots, s_n$ real and non-negative such that $\sum_{i=1}^n s_i = 1$ (convexity condition), defining $\underline{G} = \sum_{i=1}^n s_i(X_i - \underline{P}(X_i)) - (X_0 - \underline{P}(X_0))$, $\sup \underline{G} \geq 0$.*[1]

Any coherent lower prevision is convex, since Definition 3 is obtained from Definition 2 adding the constraint $\sum_{i=1}^n s_i = s_0 = 1$ (note that we would get a definition equivalent to Definition 3 requiring only $\sum_{i=1}^n s_i = s_0 > 0$). Conversely, a convex lower prevision does not even necessarily avoid sure loss:

**Proposition 1** *Let $\underline{P}$ be a convex lower prevision on $\mathcal{D}$ and let $0 \in \mathcal{D}$. Then $\underline{P}$ avoids sure loss iff $\underline{P}(0) \leq 0$.*

Convexity is characterised by a set of axioms if $\mathcal{D}$ has a special structure:

**Theorem 2** *Let $\underline{P} : \mathcal{D} \to \mathbb{R}$.*

*(a) If $\mathcal{D}$ is a* linear space *containing real constants, $\underline{P}$ is a convex lower prevision iff it satisfies the following axioms:*[2]

  *(T)* $\underline{P}(X+c) = \underline{P}(X) + c, \forall X \in \mathcal{D}, \forall c \in \mathbb{R}$ *(translation invariance)*

  *(M)* $\forall X, Y \in \mathcal{D}$, if $Y \leq X$ then $\underline{P}(Y) \leq \underline{P}(X)$ *(monotonicity)*

  *(C)* $\underline{P}(\lambda X + (1-\lambda)Y) \geq \lambda \underline{P}(X) + (1-\lambda)\underline{P}(Y), \forall X, Y \in \mathcal{D}, \forall \lambda \in [0,1]$ *(concavity).*

*(b) If $\mathcal{D}$ is a* convex cone, *$\underline{P}$ is a convex lower prevision iff it satisfies (C) and*

  *(M1)* $\forall \mu \in \mathbb{R}, \forall X, Y \in \mathcal{D}$, if $X \geq Y + \mu$ then $\underline{P}(X) \geq \underline{P}(Y) + \mu$.

**Proposition 2** *Some properties of convex lower previsions.*

---

[1] The term 'convex' in 'convex prevision' refers to the convexity condition $\sum_{i=1}^n s_i = 1$ ($s_i \geq 0$), which distinguishes convex lower (upper) previsions from coherent lower (upper) previsions (cf. Definitions 2, 3 and 7) and convex natural extensions from natural extensions (cf. Definition 4 and Section 3.2.1). The term 'convex prevision' is therefore unrelated with convexity or concavity properties of previsions as real functions.

[2] (T) and (M) can be replaced by $\underline{P}(X) - \underline{P}(Y) \leq \sup(X - Y), \forall X, Y \in \mathcal{D}$.

*(a) (Convergence theorem) Let $\{\underline{P}_j\}_{j=1}^{+\infty}$ be a sequence of lower previsions, convex on $\mathcal{D}$ and such that $\forall X \in \mathcal{D}$ there exists $\lim_{j \to +\infty} \underline{P}_j(X) = \underline{P}(X)$. Then $\underline{P}$ is convex on $\mathcal{D}$.*

*(b) (Convexity theorem) If $\underline{P}_1$ and $\underline{P}_2$ are convex lower previsions on $\mathcal{D}$, so is $\underline{P}(X) = \lambda \underline{P}_1(X) + (1-\lambda)\underline{P}_2(X)$, $\forall \lambda \in [0,1]$.*

*Let $\underline{P}$ be a convex lower prevision on $\mathcal{D}$. The following properties hold (whenever all random variables involved are in $\mathcal{D}$):*

*(c) If $\underline{P}(0) \geq 0$, $\underline{P}(\lambda X) \geq \lambda \underline{P}(X)$, $\forall \lambda \in [0,1]$ and $\underline{P}(\lambda X) \leq \lambda \underline{P}(X)$, $\forall \lambda > 1$*

*(d) $\underline{P}(0) + \inf X \leq \underline{P}(X) \leq \underline{P}(0) + \sup X$*

*(e) $\forall \mu \in \mathbb{R}$, $\underline{P}^*(X) = \underline{P}(X) + \mu$ is convex on $\mathcal{D}$.*

Properties (a) and (b), which are quite analogous to corresponding properties of coherent previsions and previsions avoiding sure loss [14], point out ways of obtaining new convex lower previsions from given ones. Property (c) shows that convexity is compatible with lack of positive homogeneity, but requires the condition $\underline{P}(0) \geq 0$. Property (d) highlights a sore point of convexity: $\underline{P}(X)$ need not belong to the closed interval $[\inf X, \sup X]$ (*internality* may fail).[3]

Property (d) suggests that internality could be restored imposing $\underline{P}(0) = 0$, if $0 \notin \mathcal{D}$; by (e), if $0 \in \mathcal{D}$ and $\underline{P}(0) \neq 0$, then $\underline{P}^*(X) = \underline{P}(X) - \underline{P}(0)$ is convex and $\underline{P}^*(0) = 0$. Requiring $\underline{P}(0) = 0$ is also the only choice to make $\underline{P}$ avoid sure loss (Proposition 1), while assuring that (c) holds.

Thinking of the meaning of a lower prevision, it appears extremely reasonable to add condition $\underline{P}(0) = 0$ to convexity: it would be at least weird to give an estimate (even imprecise) of the non-random variable 0 which is other than zero.

### 3.2 Convex Natural Extension

Before considering the stronger class of centered convex previsions, we introduce the notion of convex natural extension, which is strictly related to convexity.

**Definition 4** *Let $\underline{P} : \mathcal{D} \to \mathbb{R}$ be a lower prevision, $Z$ an arbitrary (bounded) random variable. Define $g_h = s_h(X_h - \underline{P}(X_h))$, $L = \{\alpha : Z - \alpha \geq \sum_{i=1}^n g_i$, for some $n \geq 1, X_i \in \mathcal{D}, s_i \geq 0$, with $\sum_{i=1}^n s_i = 1\}$. $\underline{E}_c(Z) = \sup L$ is termed* convex natural extension[4] *of $\underline{P}$ on $Z$.*

It is clear that $L$ is always non-empty (putting $n = 1$, $s_1 = 1$, $X_1 = X \in \mathcal{D}$ in its definition, $\alpha \in L$ for $\alpha \leq \inf Z - \sup X + \underline{P}(X)$), while $\underline{E}_c(Z)$ can in general be infinite. This situation is characterised in the following Proposition 3.

---

[3]Non-internality cannot anyway be two-sided: if there exists $X \in \mathcal{D}$ such that $\underline{P}(X) > \sup X$ ($\underline{P}(X) < \inf X$), then $\underline{P}(Y) > \inf Y$ ($\underline{P}(Y) < \sup Y$), $\forall Y \in \mathcal{D}$. This is easily seen applying Definition 3, with $n = 2$, $\{X_0, X_1\} = \{X, Y\}$.

[4]The reason why $\underline{E}_c$ is termed 'extension' appears from the later Theorem 3, especially (d).

**Definition 5** $\underline{P} : \mathcal{D} \to \mathbb{R}$ *is a lower prevision that* avoids unbounded sure loss *on $\mathcal{D}$ iff there exists $k \in \mathbb{R}$ such that, for all $n \in \mathbf{N}^+$, $\forall X_1, \ldots, X_n \in \mathcal{D}$, $\forall s_1, \ldots, s_n$ real and non-negative with $\sum_{i=1}^n s_i = 1$, defining $\underline{G} = \sum_{i=1}^n s_i(X_i - \underline{P}(X_i))$, $\sup \underline{G} \geq k$.*

**Remark 1** *Definition 5 generalises Definition 1: $\underline{P}$ avoids unbounded sure loss if and only if $\underline{P} + k$ avoids sure loss for some $k \in \mathbb{R}$, since the last inequality in Definition 5 may be written as $\sup \sum_{i=1}^n s_i(X_i - (\underline{P}(X_i) + k)) \geq 0$ and the constraint $\sum_{i=1}^n s_i = 1$ is not restrictive for Definition 1. Note also that if $\underline{P} + k$ avoids sure loss, then so does $\underline{P} + h$, $\forall h \leq k$. Therefore, when $\underline{P}$ avoids unbounded sure loss, defining $\overline{k} = \sup \{k \in \mathbb{R} : \underline{P} + k \text{ avoids sure loss}\}$, $\underline{P}$ avoids sure loss too whenever $\overline{k} \geq 0$. As a further remark, it can be seen that the constraint $\sum_{i=1}^n s_i = 1$ is essential in Definition 5: wiping it out would make Definition 5 equivalent to Definition 1.*

**Proposition 3** $\underline{E}_c(Z)$ *is finite, whatever is Z, iff $\underline{P}$ avoids unbounded sure loss.*

**Proof.** Suppose first that $\underline{P}$ avoids unbounded sure loss and for an arbitrary $Z$ let $\alpha \in L$. Then $Z - \alpha \geq \sum_{i=1}^n s_i(X_i - \underline{P}(X_i))$ for some $X_1, \ldots, X_n \in \mathcal{D}$ and $s_1, \ldots, s_n \geq 0$ with $\sum_{i=1}^n s_i = 1$, and hence $\sup Z - \alpha \geq \sup \sum_{i=1}^n s_i(X_i - \underline{P}(X_i)) \geq k$, using Definition 5 at the last inequality. Therefore $\underline{E}_c(Z) \leq \sup Z - k$.

Conversely, suppose now that $\underline{P}$ does not avoid unbounded sure loss. Therefore, for each $k \in \mathbb{R}$ there are $X_1, \ldots, X_n \in \mathcal{D}$ and $s_1, \ldots, s_n \geq 0$ with $\sum_{i=1}^n s_i = 1$ such that $\sum_{i=1}^n s_i(X_i - \underline{P}(X_i)) < k \leq Z - (-k + \inf Z)$. This implies, for any $Z$, $-k + \inf Z \in L$ and, by the arbitrariness of $k$, $\underline{E}_c(Z) = +\infty$. □

The condition of avoiding unbounded sure loss is rather mild. For instance, it clearly holds whenever $\mathcal{D}$ is finite. It is also implied by convexity, as shown by the following proposition, while the converse implication is generally not true.

**Proposition 4** *If $\underline{P} : \mathcal{D} \to \mathbb{R}$ is convex, it avoids unbounded sure loss.*

**Proof.** Choose arbitrarily $X_1, \ldots, X_n \in \mathcal{D}$ and $s_1, \ldots, s_n \geq 0$ such that $\sum_{i=1}^n s_i = 1$ in Definition 5. Given $X_0 \in \mathcal{D}$, use convexity to write $0 \leq \sup \{\sum_{i=1}^n s_i(X_i - \underline{P}(X_i)) - (X_0 - \underline{P}(X_0))\} \leq \sup \{\sum_{i=1}^n s_i(X_i - \underline{P}(X_i))\} - (\inf X_0 - \underline{P}(X_0))$, and hence $\sup \{\sum_{i=1}^n s_i(X_i - \underline{P}(X_i))\} \geq k = \inf X_0 - \underline{P}(X_0)$. □

We state now some properties of the convex natural extension. An indirect characterisation of the convex natural extension will be given in Theorem 5.

**Theorem 3** *Let $\underline{P} : \mathcal{D} \to \mathbb{R}$ be a lower prevision which avoids unbounded sure loss and $\underline{E}_c$ its convex natural extension. Then*

(a) *$\underline{E}_c$ is a convex prevision on $L$ and $\underline{E}_c(X) \geq \underline{P}(X), \forall X \in \mathcal{D}$*

(b) *$\underline{P}$ is convex if and only if $\underline{E}_c = \underline{P}$ on $\mathcal{D}$*

(c) *If $\underline{P}^*$ is a convex prevision on $L$ such that $\underline{P}^*(X) \geq \underline{P}(X)$ $\forall X \in \mathcal{D}$, then $\underline{P}^*(Z) \geq \underline{E}_c(Z), \forall Z \in L$*

(d) *If $\underline{P}$ is convex, $\underline{E}_c$ is the minimal convex extension of $\underline{P}$ to $\mathcal{L}$*

(e) *$\underline{P}$ avoids sure loss on $\mathcal{D}$ if and only if $\underline{E}_c$ avoids sure loss on $\mathcal{L}$.*

### 3.2.1 The Role of the Convex Natural Extension

The properties of $\underline{E}_c$ closely resemble those of the *natural extension $\underline{E}$* [14] of a lower prevision $\underline{P}$, whose definition differs from that of $\underline{E}_c$ only for the lack of the constraint $\sum_{i=1}^{n} s_i = 1$. In particular, as $\underline{E}$ characterises coherence of $\underline{P}$ ($\underline{P}$ is coherent iff $\underline{E}$ coincides with $\underline{P}$ on $\mathcal{D}$), $\underline{E}_c$ characterises convexity of $\underline{P}$.

Property (d) lets us extend $\underline{P}$ to *any* $\mathcal{D}' \supset \mathcal{D}$ (maintaining convexity) by considering the restriction of $\underline{E}_c$ to $\mathcal{D}'$. Moreover, (e) guarantees that $\underline{E}_c$ inherits the condition of avoiding sure loss when $\underline{P}$ satisfies it.

It is well known that the natural extension is finite iff $\underline{P}$ avoids sure loss, and when finite it can correct $\underline{P}$ into a coherent assessment in a canonical way. Analogously, the convex natural extension is finite iff $\underline{P}$ avoids unbounded sure loss, and can be used to correct $\underline{P}$ into a convex assessment, although property (e) warns us that $\underline{E}_c$ will still incur sure loss if $\underline{P}$ does so. This problem can be solved using Proposition 2, (e): $\underline{P}^*(X) = \underline{E}_c(X) - \underline{E}_c(0)$ is a correction of $\underline{P}$ which avoids sure loss by Proposition 1, as $\underline{P}^*(0) = 0$. This also means that $\underline{P}^*$ is a centered convex prevision by Definition 6 in the next section.

Alternatively, the convex natural extension may be employed to correct an assessment $\underline{P}$ which avoids unbounded sure loss (but not sure loss) into $\underline{P}'$, which avoids sure loss but is not necessarily convex. In fact, $\underline{P} + h$ avoids sure loss $\forall h \leq \overline{k} < 0$ (cf. Remark 1). Since it can be shown that $\overline{k} = -\underline{E}_c(0)$, it ensues that $\underline{E}_c(0)$ is the minimum $k$ to be subtracted from $\underline{P}$ to make $\underline{P}' = \underline{P} - k$ avoid sure loss.

Hence, the convex natural extension points out ways of correcting an assessment incurring (bounded) sure loss into one avoiding sure loss, a problem which cannot be answered using the natural extension. These corrections can be applied in several interesting situations, including, as already noted, the case of a finite $\mathcal{D}$.

## 3.3 Centered Convex Previsions and Envelope Theorems

The considerations at the end of Section 3.1 lead us naturally to the following stronger notion of centered convexity:

**Definition 6** *A lower prevision $\underline{P}$ on domain $\mathcal{D}$ ($0 \in \mathcal{D}$) is* centered convex *(C-convex, in short) iff it is convex and $\underline{P}(0) = 0$.*[5]

**Proposition 5** *Let $\underline{P}$ be a centered convex lower prevision on $\mathcal{D}$. Then*

(a) *$\underline{P}$ has a convex natural extension (hence at least one centered convex extension) on any $\mathcal{D}' \supset \mathcal{D}$*

---

[5]As shown in [12], we obtain an equivalent definition of centered convex lower prevision by requiring $\underline{P}(0) = 0$ and relaxing the convexity condition $\sum_{i=1}^{n} s_i = s_0 > 0$ to $\sum_{i=1}^{n} s_i \leq s_0$.

*(b)* $\underline{P}(\lambda X) \geq \lambda \underline{P}(X)$, $\forall \lambda \in [0,1]$, $\underline{P}(\lambda X) \leq \lambda \underline{P}(X)$, $\forall \lambda \in ]-\infty, 0[ \cup ]1, +\infty[$

*(c)* $\inf X \leq \underline{P}(X) \leq \sup X$, $\forall X \in \mathcal{D}$

*(d)* $\underline{P}$ *avoids sure loss.*

*Besides, the convergence and convexity theorems hold for C-convex previsions too (replacing 'convex' with 'centered convex' in Proposition 2, (a) and (b)).*

Properties (a)÷(d) show that centered convexity is significantly closer to coherence than convexity: C-convex lower previsions are a special class of previsions avoiding sure loss, retaining several properties of coherence and the extension property of convexity, but not requiring positive homogeneity.

A convex prevision $\underline{P}$ which is not centered may still be avoiding sure loss, if $\underline{P}(0) < 0$ (Proposition 1), but in general it is only warranted by Proposition 4 that it avoids unbounded sure loss, a very weak consistency requirement.

**Remark 2** (Convexity and n-coherence) *The consistency notion of n-coherence is discussed in [14], Appendix B, illustrating how it can be appropriate for certain 'bounded rationality' models. If the model does not require positive homogeneity, n-coherence alone is inadequate: 1-coherence is too weak, being equivalent to the internality condition (c) in Proposition 5, 2-coherence is too strong, as on linear spaces it is equivalent to two axioms, one of which is positive homogeneity [14]. As a matter of fact, C-convex previsions are a special class of 1-coherent (but not necessarily 2-coherent) previsions.*

An indirect comparison among convexity, centered convexity and coherence is given by their corresponding envelope theorems. We firstly recall that it was proved in [14] that any lower envelope of coherent lower previsions is coherent. Here is the parallel statement for convex lower previsions, while the generalisation of Theorem 1 (lower envelope theorem) comes next.

**Proposition 6** *Let $\mathcal{P}$ be a set of convex lower previsions all defined on $\mathcal{D}$. If $\underline{P}(X) = \inf_{\underline{Q} \in \mathcal{P}} \{\underline{Q}(X)\}$ is finite $\forall X \in \mathcal{D}$, $\underline{P}$ is convex on $\mathcal{D}$.*

**Theorem 4** (Generalised envelope theorem) *$\underline{P}$ is convex on $\mathcal{D}$ iff there exist a set $\mathcal{P}$ of coherent precise previsions on $\mathcal{D}$ and a function $\alpha : \mathcal{P} \to \mathbb{R}$ such that:*

*(a)* $\underline{P}(X) = \inf_{P \in \mathcal{P}} \{P(X) + \alpha(P)\}$, $\forall X \in \mathcal{D}$   (inf *is attained*).

*Moreover, $\underline{P}$ is centered convex iff ($0 \in \mathcal{D}$ and) both (a) and the following (b) hold:*

*(b)* $\inf_{P \in \mathcal{P}} \{\alpha(P)\} = 0$   (inf *is attained*).

A result similar to Theorem 4 was proved in risk measurement theory [4], requiring $\mathcal{D}$ to be a linear space. The proof of Theorem 4, given in [12] in the framework of imprecise prevision theory, is simpler and imposes no structure on $\mathcal{D}$.

**Remark 3** *In particular, the constructive implication of the theorem (for convex previsions) enables us to obtain convex previsions as lower envelopes of translated precise previsions. Its proof follows easily from Proposition 6 and Proposition 2, (e): every precise prevision P is convex and so is $P + \alpha(P)$, by Proposition 2, (e); $\inf_{P \in \mathcal{P}} \{P(X) + \alpha(P)\}$ is a convex prevision by Proposition 6.*

**Remark 4** *Let $\underline{P}$ be a lower prevision and $\mathcal{S}$ the set of all coherent precise previsions on $\mathcal{L}$. Define also $\mathcal{M}(\underline{P}) = \{(Q,r) \in \mathcal{S} \times \mathbb{R} : Q(X) + r \geq \underline{P}(X), \forall X \in \mathcal{D}\}$. It ensues from Theorem 4 that convexity of $\underline{P}$ can be equivalently characterised by the condition $\underline{P}(X) = \inf\{Q(X) + r : (Q,r) \in \mathcal{M}(\underline{P})\}$ $\forall X \in \mathcal{D}$; C-convexity can be characterised by adding the constraint $\inf\{r : \exists Q \in \mathcal{S} : (Q,r) \in \mathcal{M}(\underline{P})\} = 0$ (cf. also the following Theorem 5, where the lower envelope concerns all $X \in \mathcal{L}$).*

The envelope theorem characterisations of convexity, centered convexity and coherence differ about the role of function $\alpha$, which is unconstrained with convexity, non-negative and such that $\min \alpha = 0$ with centered convexity, identically equal to zero with coherence (in this case Theorem 4 reduces to Theorem 1).

The result in the next theorem characterises the convex natural extension as the lower envelope of a set of translated coherent precise previsions and can be proved in a way similar to the natural extension theorem in [14], Section 3.4.

**Theorem 5** *Let $\underline{P}$ be a lower prevision on $\mathcal{D}$ which avoids unbounded sure loss and define $\mathcal{S}$ and $\mathcal{M}(\underline{P})$ as in Remark 4. Then, $\mathcal{M}(\underline{P}) = \mathcal{M}(\underline{E}_c)$ and $\underline{E}_c(X) = \inf\{Q(X) + r : (Q,r) \in \mathcal{M}(\underline{P})\}, \forall X \in \mathcal{L}$.*

## 4 Some Applications

We have seen so far that convexity may help in correcting several inconsistent assessments. As noted in Section 3.2.1, its usefulness in this problem is essentially instrumental: we may easily go further and arrive at a centered convex correction, which guarantees a more satisfactory degree of consistency.

Turning to other problems, some uncertainty modelisations give rise to convex previsions, as in the examples which follow. We emphasise that we do not maintain that these models are reasonable, but simply that they are sometimes adopted in practice, and that convexity supplies a conceptual framework for them.

**Example 1** (Overly prudential assessments) *Persons or institutions which have to evaluate the random variables in a set $\mathcal{D}$ are often unfamiliar with uncertainty theories. In this case, a solution is to gather n experts and ask each of them to formulate a precise prevision (or an expectation) for all $X \in \mathcal{D}$. Choosing $\underline{P}(X) = \min_{i=1,\ldots,n} P_i(X), \forall X$ (where $P_i$ is expert i's evaluation) as one's own opinion is an already prudential way of pooling the experts' opinions, and originates a coherent lower prevision. Some more caution or lack of confidence toward some experts may lead to replacing every $P_i$ with $P_i^* = P_i - \alpha_i$ before performing the* *minimum, where $\alpha_i \geq 0$ measures in some way the final assessor's personal caution or his/her (partially) distrusting expert i. By Theorem 4, $\underline{P}^* = \min_{i=1,\ldots,n} \underline{P}_i^*$ is convex (cf. Remark 3). More generally, $\underline{P}^*$ is of course convex also when the sign of the $\alpha_i$ is unconstrained ($\alpha_i < 0$ if, for instance, expert i's opinion is believed to be biased and below the 'real' prevision). It is interesting to observe that if $\alpha_i \geq 0$ for at least one i, $\underline{P}^*$ avoids sure loss too (since then $\underline{E}_c(0) \leq 0$ by Theorem 5, hence $\underline{E}_c$ avoids sure loss by Proposition 1, and so does $\underline{P}^*$ by Theorem 3, (e)). In particular, the following situation may be not unusual with an unexperienced assessor: $\alpha_i > 0$ for some i, and $0 \notin \mathcal{D}$, because the assessor thinks that no expert is needed to evaluate 0, he himself can assign, of course, $\underline{P}^*(0) = 0$. If such is the case, the extension of $\underline{P}^*$ on $\mathcal{D} \cup \{0\}$ keeps on avoiding sure loss, as is easily seen, but is generally not convex (to see this with a simple example, suppose $X \in \mathcal{D}$, $\underline{P}^*(X) < \inf X$ and use the result in footnote 3 to obtain that $\underline{P}^*(0) < 0$ is then necessary for convexity).*

In the following example and in Section 4.1 we shall refer to upper previsions, to which the theory developed so far extends with mirror-image modifications. We report the conjugates of Definition 3 and Theorem 4.

**Definition 7** $\overline{P} : \mathcal{D} \to \mathbb{R}$ is a convex upper prevision *on $\mathcal{D}$ iff, for all $n \in \mathbf{N}^+$, $\forall X_0, X_1, \ldots, X_n \in \mathcal{D}$, $\forall s_1, \ldots, s_n$ real and non-negative such that $\sum_{i=1}^n s_i = 1$ (convexity condition), defining $\overline{G} = \sum_{i=1}^n s_i(\overline{P}(X_i) - X_i) - (\overline{P}(X_0) - X_0)$, $\sup \overline{G} \geq 0$.*

**Theorem 6** *$\overline{P}$ is convex on its domain $\mathcal{D}$ iff there exist a set $\mathcal{P}$ of coherent precise previsions (all defined on $\mathcal{D}$) and a function $\alpha : \mathcal{P} \to \mathbb{R}$ such that:*

*(a) $\overline{P}(X) = \sup_{P \in \mathcal{P}} \{P(X) + \alpha(P)\}$, $\forall X \in \mathcal{D}$ (sup is attained).*

*Moreover, $\overline{P}$ is centered convex iff ($0 \in \mathcal{D}$ and) both (a) and the following (b) hold:*

*(b) $\sup_{P \in \mathcal{P}} \{\alpha(P)\} = 0$ (sup is attained).*

**Example 2** (Supremum preserving functions) *Let $\mathbb{P} = \{\omega_i\}_{i \in I}$ be a (not necessarily finite) set of exhaustive non-impossible elementary events or* atoms, *i.e. $\omega_i \neq \varnothing \; \forall i \in I$, $\cup_{i \in I} \omega_i = \Omega$, $\omega_i \cap \omega_j = \varnothing$ if $i \neq j$. Given a function $\pi : \mathbb{P} \to [0,1]$, define $\Pi : 2^{\mathbb{P}} - \{\varnothing\} \to [0,1]$ ($2^{\mathbb{P}}$ is the powerset of $\mathbb{P}$) by*

$$\Pi(A) = \sup_{\omega_i \in A} \{\pi(\omega_i)\}, \forall A \in 2^{\mathbb{P}} - \{\varnothing\}. \tag{1}$$

*As well-known, if $\pi$ is normalised (i.e., $\sup \pi = 1$) and extended to $\varnothing$ putting $\pi(\varnothing)(= \Pi(\varnothing)) = 0$, $\Pi$ is a normalised possibility measure, a special case of coherent upper probability [3]. Without these additional assumptions, $\Pi$ is a convex upper probability. To see this, define for $i \in I$, $P_i(\omega_i) = 1$, $P_i(\omega_j) = 0 \; \forall j \neq i$, $\alpha_i = \pi(\omega_i) - 1$, and extend (trivially) each $P_i$ to $2^{\mathbb{P}}$. It is not difficult to see that $\Pi(A) = \sup_{i \in I} \{P_i(A) + \alpha_i\}$, $\forall A \in 2^{\mathbb{P}}$ and therefore $\Pi$ is convex by Theorem 6. If*

$\sup \pi < 1$, $\Pi$ *has the unpleasant property that* $\Pi(\Omega) < 1$, *and also* $\Pi(\varnothing) < 0$ *(this means that* $\Pi$ *incurs sure loss and is not C-convex). Functions similar to these kinds of unnormalised possibilities were considered in the literature relating possibility and fuzzy set theories, and their unsatisfactory properties were already pointed out (see e.g. [9], Section 2.6 and the references quoted therein).*

## 4.1 Convex Risk Measures

Further applications of convex imprecise previsions are suggested by the fact that they do not necessarily require positive homogeneity, as appears from Proposition 5, (b). Considering the well-known behavioural interpretation of lower (and upper) previsions [14], it is intuitively clear that applications could be generally related to situations of risk aversion, because of which an agent's supremum buying price for the random quantity $\lambda X$ might be less than $\lambda$ times his/her supremum buying price for $X$, when $\lambda > 1$.

In this section we shall discuss an application to (financial) risk measurement. The literature on risk measures is quite large, as this topic is very important in many financial, banking or insurance applications. Formally, a risk measure is a mapping $\rho$ from a set $\mathcal{D}$ of random variables into $\mathbb{R}$. Therefore $\rho$ associates a real number $\rho(X)$ to every $X \in \mathcal{D}$, which should determine how 'risky' $X$ is, and whether it is acceptable to buy or hold $X$. Intuitively, $X$ should be acceptable (not acceptable) if $\rho(X) \le 0$ (if $\rho(X) > 0$), and $\rho(X)$ should determine the maximum amount of money which could be subtracted from $X$, keeping it acceptable (the minimum amount of money to be added to $X$ to make it acceptable).

Traditional risk measures, like Value-at-Risk (*VaR*) – probably the most widespread – require assessing (at least) a distribution function for each $X \in \mathcal{D}$; often, a joint normal distribution is assumed [8]. Quite recently, other risk measures were introduced, which do not require assessing exactly one precise probability distribution for each $X \in \mathcal{D}$, and are therefore appropriate also in situations where conflicting or insufficient information is available. Precisely, coherent risk measures were defined in a series of papers (including [1, 2]) using a set of axioms (among these positive homogeneity), and assuming that $\mathcal{D}$ is a linear space. In these papers, coherent risk measures were not related with imprecise previsions theory, while this was done in [11, 13]; see also [10] for a general approach to these and other theories. Convex risk measures were introduced in [4, 5, 7] as a generalisation of coherent risk measures which does not require the positive homogeneity axiom. We report the definition in [5]:

**Definition 8** *Let $\mathcal{V}$ be a linear space of random variables which contains real constants.* $\rho : \mathcal{V} \to \mathbb{R}$ *is a* convex risk measure *iff it satisfies the following axioms:*

*(T1)* $\forall X \in \mathcal{V}, \forall \alpha \in \mathbb{R}, \rho(X + \alpha) = \rho(X) - \alpha$ *(translation invariance)*

*(M2)* $\forall X, Y \in \mathcal{V}$, *if* $X \le Y$ *then* $\rho(Y) \le \rho(X)$ *(monotonicity)*

*(C1)* $\rho(\lambda X + (1-\lambda)Y) \leq \lambda\rho(X) + (1-\lambda)\rho(Y) \; \forall X,Y \in \mathcal{V}, \lambda \in [0,1]$ *(convexity)*.

Convex risk measures are also discussed in [6] and their potential capability of capturing risk aversion is pointed out in [5]. In a risk measurement environment, a motivation for not assuming positive homogeneity is that $\rho(\lambda X)$ may be larger than $\lambda\rho(X)$ for $\lambda > 1$ also because of *liquidity risks*: if we were to sell immediately a large amount $\lambda X$ of a financial investment, we might be forced to accept a smaller reward than $\lambda$ times the current selling price for $X$.

It was shown in [11] that risk measures can be encompassed into the theory of imprecise previsions, because a risk measure for $X$ can be interpreted as an upper prevision for $-X$:[6]

$$\rho(X) = \overline{P}(-X). \tag{2}$$

This fact was used in [11, 13] to generalise the notion of coherent risk measures to an arbitrary domain $\mathcal{D}$. An analogue generalisation can be done for convex risk measures [12], as we shall now illustrate.

**Definition 9** $\rho : \mathcal{D} \to \mathbb{R}$ *is a* convex risk measure *on $\mathcal{D}$ if and only if for all $n \in \mathbb{N}^+$, $\forall X_0, X_1, \ldots, X_n \in \mathcal{D}$, $\forall s_1, \ldots, s_n$ real and non-negative such that $\sum_{i=1}^n s_i = 1$, defining $\overline{G} = \sum_{i=1}^n s_i(X_i + \rho(X_i)) - (X_0 + \rho(X_0))$, $\sup \overline{G} \geq 0$.*

Note that Definition 9 may be obtained from Definition 7 referring to $-X$ rather than $X$, for all $X \in \mathcal{D}$.

If $\mathcal{D}$ is a linear space containing real constants, the notion in Definition 9 reduces to that in [4, 5], by the next theorem (cf. also Theorem 2, (a)):

**Theorem 7** *Let $\mathcal{V}$ be a linear space of bounded random variables containing real constants. A mapping $\rho$ from $\mathcal{V}$ into $\mathbb{R}$ is a convex risk measure according to Definition 9 iff it is a convex risk measures according to Definition 8.*

Definition 9 applies to any set $\mathcal{D}$ of random variables, unlike Definition 8, which, if $\mathcal{D}$ is arbitrary, requires embedding it in a larger linear space.

Results specular to those presented in Section 3 apply to convex risk measures. In particular, the convergence and convexity theorems (Proposition 2, (a) and (b)) hold; convex risk measures can be extended on any $\mathcal{D}' \supset \mathcal{D}$, preserving convexity; they avoid sure loss iff $\rho(0) \geq 0$ (we say that $\rho$ avoids sure loss on $\mathcal{D}$ iff $\overline{P}(-X) = \rho(X)$ avoids sure loss on $\mathcal{D}^- = \{-X : X \in \mathcal{D}\}$).

Like the general case in Section 3, it appears quite appropriate to put $\rho(0) = 0$, and hence to use *centered convex* risk measures: 0 is the unquestionably reasonable selling or buying price for $X = 0$.

**Definition 10** *A mapping $\rho$ from $\mathcal{D}$ ($0 \in \mathcal{D}$) into $\mathbb{R}$ is a* centered convex risk measure *on $\mathcal{D}$ iff $\rho$ is convex and $\rho(0) = 0$.*

---

[6] We assume that the time gap between the buying and selling time of $X$ is negligible (if not, we should introduce a discounting factor in (2)). This simplifies the sequel, without substantially altering the conclusions.

Centered convex risk measures have further nice additional properties, corresponding to those of centered convex lower previsions: they always avoid sure loss, and are such that $-\sup X \leq \rho(X) \leq -\inf X, \forall X \in \mathcal{D}$.

This condition corresponds to internality ((c) of Proposition 5), and is a rationality requirement for risk measures: for instance, $\rho(X) > -\inf X$ would mean that to make $X$ acceptable we require adding to it a sure number ($\rho(X)$) higher than the maximum loss $X$ may cause.

A centered convex risk measures $\rho$ is not necessarily positively homogeneous:

$$\rho(\lambda X) \geq \lambda \rho(X), \forall \lambda \geq 1. \tag{3}$$

A notion of convex natural extension may also be given for centered convex (or convex) risk measures and its properties correspond to those listed in Theorem 3. When finite, it gives in particular a standard way of 'correcting' other kinds of risk measures into convex risk measures.[7]

The generalised envelope theorem is obtained from the statement of Theorem 6 replacing $\overline{P}(X)$ and $\underline{P}(X)$ with, respectively, $\rho(X)$ and $P(-X)$.

Examples of convex risk measures may be found in [4, 5, 12].

## 5 Conclusions

In this paper we studied convex and centered convex previsions in the framework of Walley's theory of imprecise previsions. Convex previsions do not necessarily satisfy minimal consistency requirements, but are useful in generalising natural extension-like methods of correcting inconsistent assessments and in providing a conceptual framework for some uncertainty models. Centered convex previsions are in a sense intermediate between avoiding sure loss and coherence: their properties are closer to coherence than those of a generic prevision that avoids sure loss, but are also compatible with lack of positive homogeneity. Because of this, they are potentially useful at least in models which incorporate some forms of risk aversion. We outlined a risk measurement application, where they lead to defining convex risk measures, and believe that several applications of convex imprecise previsions are still to be explored. It might also be interesting to investigate if and how convex previsions can be generalised in a conditional environment, or when allowing unbounded random variables.

---

[7]Note that this is always possibile if $\mathcal{D}$ is finite (cf. Section 3.2.1).

**Renato Pelessoni** and **Paolo Vicig** are with the Department of Applied Mathematics 'B. de Finetti', University of Trieste, Piazzale Europa 1, I-34127 Trieste, Italy. E-mails: renato.pelessoni@econ.units.it, paolo.vicig@econ.units.it